\documentclass[11pt,twoside,leqno]{aomamlt2e} 
  \pageno{941}
\received{May 22, 2004}

\numberwithin{equation}{section}

\input xy
\xyoption{all}

\DeclareMathOperator{\Aut}{Aut}

\newcommand{\fhi}{\varphi}
\newcommand{\g}{\mathfrak{g}}
\newcommand{\f}{\mathfrak{f}}
\newcommand{\p}{\mathfrak{p}}
\newcommand{\OO}{\mathcal{O}}
\newcommand{\FF}{\mathcal{F}}
\newcommand{\R}{\mathbb{R}}
\newcommand{\C}{\mathbb{C}}
\DeclareMathOperator{\Iso}{Iso}
\DeclareMathOperator{\Gl}{Gl}
\DeclareMathOperator{\rank}{rank}

\newtheorem{thm}{Theorem}[section]

\newtheorem{lemma}[thm]{Lemma}
\newtheorem{cor}[thm]{Corollary}
\newtheorem{prop}[thm]{Proposition}

 \newcounter{notes}

\begin{document}
\currannalsline{164}{2006}

\title{Isometric actions of simple Lie groups\\ on pseudoRiemannian manifolds}

\shorttitle{Isometric actions of simple Lie groups}

\acknowledgement{Research supported by SNI-M\'exico and CONACYT Grant 44620.}
\author{Raul Quiroga-Barranco}
\institution{Centro de Investigaci\'on en Matem\'aticas, Guanajuato, Gto., M\'exico
\\
\email{quiroga@cimat.mx
}}

\centerline{\bf Abstract}
\vglue6pt

Let $M$ be a connected compact pseudoRiemannian manifold acted
upon topologically transitively and isometrically by a connected
noncompact simple Lie group $G$. If $m_0, n_0$ are the dimensions
of the maximal lightlike subspaces tangent to $M$ and $G$,
respectively, where $G$ carries any bi-invariant metric, then we
have $n_0 \leq m_0$. We study $G$-actions that satisfy the
condition $n_0 = m_0$. With no rank restrictions on $G$, we prove
that $M$ has a finite covering $\widehat{M}$ to which the
$G$-action lifts so that $\widehat{M}$ is $G$-equivariantly
diffeomorphic to an action on a double coset $K\backslash
L/\Gamma$, as considered in Zimmer's program, with $G$ normal in
$L$ (Theorem A). If $G$ has finite center and $\rank_\R(G)\geq 2$,
then we prove that we can choose $\widehat{M}$ for which $L$ is
semisimple and $\Gamma$ is an irreducible lattice (Theorem B). We
also prove that our condition $n_0 = m_0$ completely
characterizes, up to a finite covering, such double coset
$G$-actions (Theorem~C). This describes a large family of double
coset $G$-actions and provides a partial positive answer to the
conjecture proposed in Zimmer's program.

\vglue-23pt\phantom{up}

\section{Introduction}
\label{section-intro}
\vglue-6pt

In this work, $G$ will denote a connected noncompact simple Lie
group and $M$ a connected   smooth manifold, which is assumed to be compact unless otherwise stated.
Moreover, we will assume that $G$ acts smoothly, faithfully and preserving a
finite measure on $M$. We will assume that these conditions are
satisfied unless stated otherwise. There are several known
examples of such actions that also preserve some geometric
structure and all of them are essentially of an algebraic nature
(see \cite{Zimmer-prog} and \cite{FK-survey}). Some of such
examples are constructed from homomorphisms $G\hookrightarrow L$
into Lie groups $L$ that admit a (cocompact) lattice $\Gamma$. For
such setup, the $G$-action is then the one by left translations on
$K\backslash L/\Gamma$, where $K$ is some compact subgroup of
$C_L(G)$. Moreover, if $L$ is semisimple and $\Gamma$ is
irreducible, then the $G$-action is ergodic. This family of
examples is a fundamental part in the questions involved in
studying and classifying $G$-actions. In his program to study such
actions, Robert Zimmer has proposed the problem of determining to
what extent a general $G$-action on $M$ as above is (or at least
can be obtained from) an algebraic action, which includes the
examples $K\backslash L/\Gamma$ as above (see \cite{Zimmer-prog}).

Our goal is to make a contribution to Zimmer's program within the
context of pseudoRiemannian geometry. Hence, from now on, we
consider $M$ furnished with a smooth pseudoRiemannian metric and
assume that $G$ acts by isometries of the metric. Note that $G$
also preserves the pseudoRiemannian volume on~$M$, which is finite
since $M$ is compact.

One of the first things we want to emphasize is the fact that $G$
itself can be naturally considered as a pseudoRiemannian manifold.
In fact, $G$ admits bi-invariant pseudoRiemannian metrics and all
of them can be described in terms of the Killing form (see
\cite{Paco-thesis} and \cite{Bader-Nevo}). So it is natural to
inquire about the relationship of the pseudoRiemannian invariants
of both $G$ and $M$. The simplest one to consider is the
signature, which from now on we will denote with $(m_1,m_2)$ and
$(n_1,n_2)$ for $M$ and $G$, respectively, where we have chosen
some bi-invariant pseudoRiemannian metric on $G$. Our notation is
such that the first number corresponds to the dimension of the
maximal timelike tangent subspaces and the second number to the
dimension of the maximal spacelike tangent subspaces. We will also
denote $m_0 = \min(m_1,m_2)$ and $n_0 = \min(n_1,n_2)$, which are
the dimensions of maximal lightlike tangent subspaces for $M$ and
$G$, respectively. We observe that the signature $(n_1,n_2)$
depends on the choice of the metric on $G$. However, as it was
remarked by Gromov in \cite{Gromov-rigid}, if $(n_1,n_2)$
corresponds to the metric given by the Killing form, then any
other bi-invariant pseudoRiemannian metric on $G$ has signature
given by either $(n_1,n_2)$ or $(n_2,n_1)$. In particular, $n_0$
does not depend on the choice of the bi-invariant metric on $G$,
so it only depends on $G$ itself. For these numbers, it is easy to
check the following inequality. A proof is given later on in
Lemma~\ref{signatures2}.

\vglue-24pt
\phantom{up}

\begin{lemma}
\label{signatures} For $G$ and $M$ as before{\rm ,} we have $n_0 \leq
m_0$.
\end{lemma}
\phantom{up}
\vglue-24pt

The goal of this paper is to obtain a complete description, in
algebraic terms, of the manifolds $M$ and the $G$-actions that
occur when the equality $n_0 = m_0$ is satisfied. We will prove
the following result. We refer to \cite{Zimmer-engagement} for the
definition of engagement.

\vskip8pt {\scshape Theorem A.} 
{\it Let $G$ be a connected noncompact simple Lie group. If $G$ acts
faithfully and topologically transitively on a compact manifold
$M$ preserving a pseudoRiemannian metric such that $n_0 = m_0${\rm ,}
then the $G$-action on $M$ is ergodic and engaging{\rm ,} and there
exist\/{\rm :}\/
\vskip3pt 
\noindent {\rm (1)} a finite covering $\widehat{M}\rightarrow M${\rm ,}
\vskip3pt 
\noindent {\rm (2)}  a connected Lie group $L$ that contains $G$ as a factor{\rm ,}
\vskip3pt

\noindent {\rm (3)}  a cocompact discrete subgroup $\Gamma$ of $L$ and
    a compact subgroup $K$ of $C_L(G)${\rm ,}
\vskip3pt 
\noindent 
for which the $G$-action on $M$ lifts to $\widehat{M}$ so that
$\widehat{M}$ is $G$-equivariantly diffeomorphic to $K\backslash
L/\Gamma$. Furthermore{\rm ,} there is an ergodic and engaging
$G$-invariant finite smooth measure on $L/\Gamma$.}
\Enddemo

In other words, if the (pseudoRiemannian) geometries of $G$ and
$M$ are closely related, in the sense of satisfying $n_0 = m_0$,
then, up to a finite covering, the $G$-action is given by the
algebraic examples considered in Zimmer's program. This result
does not require any conditions on the center or real rank of $G$.

On the other hand, it is of great interest to determine the
structure of the Lie group $L$ that appears in Theorem A. For
example, one might expect to able to prove that $L$ is semisimple
and $\Gamma$ is an irreducible lattice. By imposing some
restrictions on the group $G$, in the following result we prove
that such conclusions can be obtained. In this work we adopt the
definition of irreducible lattice found in \cite{Morris}, which
applies for connected semisimple Lie groups with finite center,
even if such groups admit compact factors. We also recall that a
semisimple Lie group $L$ is called isotypic if its Lie algebra
$\mathfrak{l}$ satisfies $\mathfrak{l}\otimes \C = \mathfrak{d}
\oplus \dots \oplus \mathfrak{d}$ for some complex simple Lie
algebra $\mathfrak{d}$.

\demo{\scshape Theorem B}
{\it Let $G$ be a connected noncompact simple Lie group with finite
center and $\rank_\R(G) \geq 2$. If $G$ acts faithfully and
topologically transitively on a compact manifold $M$ preserving a
pseudoRiemannian metric such that $n_0 = m_0${\rm ,} then there exist\/{\rm :}\/
\begin{itemize}
    \item[{\rm (1)}] a finite covering $\widehat{M}\rightarrow M${\rm ,}
    \item[{\rm (2)}] a connected isotypic semisimple Lie group $L$ with
    finite center that contains $G$ as a factor{\rm ,}
    \item[{\rm (3)}] a cocompact irreducible lattice $\Gamma$ of $L$ and
    a compact subgroup $K$ of $C_L(G)${\rm ,}
\end{itemize}
for which the $G$-action on $M$ lifts to $\widehat{M}$ so that
$\widehat{M}$ is $G$-equivariantly diffeo\-morphic to $K\backslash
L/\Gamma$. Hence{\rm ,} up to fibrations with compact fibers{\rm ,} $M$ is 
$G$-equi\-variantly diffeomorphic to $K\backslash L/\Gamma$ and
$L/\Gamma$.}
\Enddemo

To better understand these results, one can look at the geometric
features of the known algebraic actions of simple Lie groups. This
is important for two reasons. To verify that there actually exist
examples of topologically transitive actions that satisfy our
condition $n_0 = m_0$, and to understand to what extent Theorems A
and B describe such examples.

First recall that every semisimple Lie group with finite center
admits cocompact lattices. However, not every such group admits an
{\it irreducible} cocompact lattice, which is a condition
generally needed to provide ergodic actions. In the work of
\cite{Johnson} one can find a complete characterization of the
semisimple groups with finite center and without compact factors
that admit irreducible lattices. Also, in \cite{Morris}, one can
find conditions for the existence of irreducible lattices on
semisimple Lie groups with finite center that may admit compact
factors. Based on the results in \cite{Johnson} and \cite{Morris}
we state the following proposition that provides a variety of
examples of ergodic pseudoRiemannian metric preserving actions for
which $n_0 = m_0$. Its proof is an easy consequence of
\cite{Johnson} and \cite{Morris}, and the remarks that follow the
statement.

\begin{prop}
\label{lattice-isotypic} Suppose that $G$ has finite center and
$\rank_\R(G)\geq 2$. Let $L$ be a semisimple Lie group with finite
center that contains $G$ as a normal subgroup. If $L$ is isotypic{\rm ,}
then $L$ admits a cocompact irreducible lattice. Hence{\rm ,} for any
choices of a cocompact irreducible lattice $\Gamma$ in $L$ and a
compact subgroup $K$ of $C_L(G)${\rm ,} $G$ acts ergodically{\rm ,} and hence
topologically transitively{\rm ,} on $K\backslash L/\Gamma$ preserving a
pseudoRiemannian metric for which $n_0 = m_0$.
\end{prop}

For the existence of the metric, we observe that there is an
isogeny between $L$ and $G\times H$ for some connected semisimple
group $H$. On a product $G\times H$, we have $K\subset HZ(G)$ and
we can build the metric from the Killing form of $\g$ and a
Riemannian metric on $H$ which is $K$-invariant on the left and
$H$-invariant on the right. For general $L$ a similar idea can be
applied.

Hence, Proposition \ref{lattice-isotypic} ensures that topological
transitivity and the condition $n_0 = m_0$, assumed by Theorems A
and B, are satisfied by a large and important family of examples,
those built out of isotypic semisimple Lie groups containing $G$ as
a normal subgroup.

A natural problem is to determine to what extent topological
transitivity and the condition $n_0 = m_0$ characterize the
examples given in Proposition \ref{lattice-isotypic}. We obtain
such a characterization in the following result.

\demo{\scshape Theorem C}
{\it Let $G$ be a connected noncompact simple Lie group with finite
center and $\rank_\R(G) \geq 2$. Assume that $G$ acts faithfully
on a compact manifold $X$. Then the following conditions are
equivalent.
\begin{itemize}
    \item[{\rm (1)}] There is a finite covering $\widehat{X} \rightarrow
    X$ for which the $G$-action on $X$ lifts to a topologically
    transitive $G$-action on $\widehat{X}$ that preserves a
    pseudoRiemannian metric satisfying $n_0 = m_0$.
    \item[{\rm (2)}] There is a connected isotypic semisimple Lie group
    $L$ with finite center that contains $G$ as a factor{\rm ,} a cocompact
    irreducible lattice $\Gamma$ of $L$ and a compact subgroup
    $K$ of $C_L(G)$ such that $K\backslash L/\Gamma$ is a finite
    covering of $X$ with $G$-equivariant covering map.
\end{itemize}}

In words, up to finite covering maps, for topologically transitive
$G$-actions on compact manifolds, to preserve a pseudoRiemannian
metric with $n_0 = m_0$ is a condition that characterizes those
algebraic actions considered in Zimmer's program corresponding to
the double cosets $K\backslash L /\Gamma$ described in (2).

In the theorems stated above we are assuming the pseudoRiemannian
manifold acted upon by $G$ to be compact. However, it is possible to
extend our arguments to finite volume manifolds if we consider
complete pseudoRiemannian structures. In Section 8 we present the
corresponding generalizations of Theorems A, B, and C that can be
thus obtained.

With the results discussed so far, we completely describe (up to
finite coverings) the isometric actions of noncompact simple Lie
groups that satisfy our geometric condition $n_0 = m_0$. Moreover,
we have actually shown that the collection of manifolds defined by
such condition is (up to finite coverings) a very specific and
important family of the examples considered in Zimmer's program:
those given by groups containing $G$ as a normal subgroup.

Given the previous remarks, we can say that we have fully described
and classified a distinguished family of $G$-actions. Nevertheless,
it is still of interest to conclude (from our classification)
results that allow us to better understand the topological and
geometric restrictions satisfied by the family of $G$-actions under
consideration. This also allows us to make a comparison with results
obtained in other works (see, for example, \cite{FK-survey},
\cite{Lubotzky-Zimmer2}, \cite{Spatzier-Zimmer},
\cite{Zimmer-full-entropy} and \cite{Zimmer-prog}). With this
respect, in the theorems below, and under our standing condition
$n_0 = m_0$, we find improvements and/or variations of important
results concerning volume preserving $G$-actions. Based on this, we
propose the problem of extending such theorems to volume preserving
$G$-actions more general than those considered here.

In the remaining of this section, we will assume that $G$ is a
connected noncompact simple Lie group acting smoothly, faithfully
and topologically transitively on a manifold $M$ and preserving a
pseudoRiemannian metric such that $n_0 = m_0$. We also assume that
either $M$ is compact or its metric is complete with finite volume.
The results stated below basically follow from Theorems A, B and C
(and their extensions to finite volume complete manifolds); the
corresponding proofs can be found in Section \ref{section-further}.

The next result is similar in spirit to Theorem A in
\cite{Spatzier-Zimmer}, but requires no rank restriction on $G$.

\begin{thm}
\label{Riemannian-split} If the $G$-action is not transitive, then
$M$ has a finite covering space $M_1$ that admits a Riemannian
metric whose universal covering splits isometrically. In particular{\rm ,}
for such metric, $M_1$ has some zeros for its sectional curvature.
\end{thm}

Observe that any algebraic $G$-action of the form $K\backslash
L/\Gamma$, as in Zimmer's program, is easily seen to satisfy the
conclusion of Theorem \ref{Riemannian-split} by just requiring $L$
to have at least two noncompact factors. Hence, one may propose the
problem of finding a condition, either geometric or dynamical, that
characterizes the conclusion of Theorem \ref{Riemannian-split} or an
analogous property.

The following result can be considered as an improved version of
Gromov's representation theorem. In this case we require a rank
restriction.

\begin{thm}
\label{gromov-rep2} Suppose $G$ has finite center and $\rank_\R(G)
\geq 2$. Then there exist a finite index subgroup $\Lambda$ of
$\pi_1(M)$ and a linear representation $\rho : \Lambda \rightarrow
\Gl(p,\R)$ such that the Zariski closure
$\overline{\rho(\Lambda)}^Z$ is a semisimple Lie group with finite
center in which $\rho(\Lambda)$ is a lattice and that contains a
closed subgroup locally isomorphic to $G$. Moreover{\rm ,} if $M$ is not
compact, then $\overline{\rho(\Lambda)}^Z$ has no compact factors.
\end{thm}

Again, we observe that all algebraic $G$-actions in Zimmer's
program, i.e. of the form $K\backslash L/\Gamma$ described before,
are easily seen to satisfy the conclusions of Theorem
\ref{gromov-rep2}. Actually, our proof depends on the fact that our
condition $n_0 = m_0$ ensures that such a double coset appears.
Still we may propose the problem of finding other conditions that
can be used to prove this more general Gromov's representation
theorem. Such a result, in a more general case, would provide a
natural semisimple Lie group in which to embed $G$ to prove that a
given $G$-action is of the type considered in Zimmer's program.

Zimmer has proved in \cite{Zimmer-full-entropy} that when
$\rank_\R(G)\geq 2$ any analytic engaging $G$-action on a manifold
$X$ preserving a unimodular rigid geometric structure has a fully
entropic virtual arithmetic quotient (see \cite{Lubotzky-Zimmer1},
\cite{Lubotzky-Zimmer2} and \cite{Zimmer-full-entropy} for the
definitions and precise statements). The following result, with our
standing assumption $n_0 = m_0$, has a much stronger conclusion than
that of the main result in \cite{Zimmer-full-entropy}. Note that a
sufficiently strong generalization of the next theorem for general
finite volume preserving actions would mean a complete solution to
Zimmer's program for finite measure preserving $G$-actions, even at
the level of the smooth category.

\begin{thm}
\label{virtual-arithmetic-quotient} Suppose $G$ and $M$ satisfy the
hypotheses of either Theorem {\rm B} or Theorem {\rm B}$'$ \/{\rm (}\/see {\rm\S
 \ref{section-further}).} Then the $G$-action on $M$ has finite
entropy. Moreover{\rm ,} there is a manifold $\widehat{M}$ acted upon by
$G$ and $G$-equivariant finite covering maps $\widehat{M}
\rightarrow A(M)$ and $\widehat{M} \rightarrow M${\rm ,} where $A(M)$ is
some realization of the maximal virtual arithmetic quotient of $M$.
\end{thm}

The organization of the article is as follows. The proof of Theorem
A relies on studying the pseudoRiemannian geometry of $G$ and $M$.
In that sense, the fundamental tools for the proof of Theorem~A are
developed in Sections \ref{section-split} and~\ref{section-G}. In
Section \ref{section-thmA} the proof of Theorem A is completed based
on the results proved up to that point and a study of a transverse
Riemannian structure associated to the $G$-orbits. The proofs of
Theorems B and C (\S\S \ref{section-thmB} and
\ref{section-thmC}) are based on Theorem A, but also rely on the
results of \cite{Stuck-Zimmer} and \cite{Zimmer-Lie}. In Section
\ref{section-further} we show how to extend Theorems A, B and C to
finite volume manifolds if we assume completeness of the
pseudoRiemannian structure involved. Section \ref{section-further}
also contains the complete proofs of Theorems
\ref{Riemannian-split}, \ref{gromov-rep2} and
\ref{virtual-arithmetic-quotient}.

I would like to thank Jes\'us \'Alvarez-L\'opez, Alberto Candel
and Dave\break Morris for useful comments that allowed to simplify the
exposition of this work.

\section{Some preliminaries on homogeneous spaces}
\label{section-prelim} 

We will need the following easy to prove
result.

\begin{lemma}
\label{finite-components} Let $H$ be a Lie group acting smoothly
and transitively on a connected manifold $X$. If for some $x_0\in
X$ the isotropy group $H_{x_0}$ has finitely many components, then
$H$ has finitely many components as well.
\end{lemma}
\Proof Let $H_{x_0} = K_0\cup\dots\cup K_r$ be the component
decomposition of $H_{x_0}$. Choose an element $k_i\in K_i$, for
every $i=0,\dots,r$.

For any given $h\in H$, let $\widehat{h}\in H_0$ be such that
$h(x_0) = \widehat{h}(x_0)$ (see \cite{Helgason}). Hence,
$\widehat{h}^{-1}h \in H_{x_0}$, so there exists $i_0$ such that
$\widehat{h}^{-1}h \in K_{i_0}$. If $\gamma$ is a continuous path
from $k_{i_0}$ to $\widehat{h}^{-1}h$, then $\widehat{h}\gamma$ is a
continuous path from $\widehat{h}k_{i_0}$ to $h$. This shows that $H
= H_0 k_0 \cup \dots H_0 k_r$. \Endproof

As an immediate consequence we obtain the following.

\begin{cor}
\label{isoX} If $X$ is a connected homogeneous Riemannian
manifold{\rm ,} then the group of isometries $\Iso(X)$ has finitely many
components. Moreover{\rm ,} the same property holds for any closed
subgroup of $\Iso(X)$ that acts transitively on $X$.
\end{cor}

The following result is a well known easy consequence of Singer's
Theorem (see \cite{Singer}). Nevertheless, we state it here for
reference and briefly explain its proof, from the results of
\cite{Singer}, for the sake of completeness.

\begin{thm}[Singer]
\label{singer-thm} Let $X$ be a smooth simply connected complete
Riemannian manifold. If the pseudogroup of local isometries has a
dense orbit{\rm ,} then $X$ is a homogeneous Riemannian manifold.
\end{thm}

\Proof By the main theorem in \cite{Singer}, we need to show that
$X$ is infinitesimally homogeneous as considered in \cite{Singer}.
The latter is defined by the existence of an isometry $A : T_x X
\rightarrow T_y X$, for any two given points $x,y \in X$, so that
$A$ transforms the curvature and its covariant derivatives (up to a
fixed order) at $x$ into those at $y$. Under our assumptions, this
condition is satisfied only on a dense subset $S$ of $X$. However,
for an arbitrary $y\in X$, we can choose $x\in S$, a sequence
$(x_n)_n \subset S$ that converges to $y$ and a sequence of maps
$A_n : T_x X \rightarrow T_{x_n} X$ that satisfy the infinitesimal
homogeneity condition. By introducing local coordinates at $x$ and
$y$, we can consider that (for $n$ large enough) the sequence
$(A_n)_n$ lies in a compact group and thus has a subsequence that
converges to some map $A : T_x X \rightarrow T_y X$. By the
continuity of the identities that define infinitesimal homogeneity
in \cite{Singer}, it is easy to show that $A$ satisfies such
identities. This proves infinitesimal homogeneity of $X$, and so $X$
is homogeneous. \hfill\qed

\section{Isometric splitting of a covering of $M$}
\label{section-split}

 We start by describing some geometric
properties of the $G$-orbits on $M$ when the condition $n_0 = m_0$
is satisfied.

\begin{prop}
\label{local-free} Suppose $G$ acts topologically transitively on
$M$ preserving its pseudoRiemannian metric and satisfying $n_0 =
m_0$. Then $G$ acts everywhere locally freely with nondegenerate
orbits. Moreover{\rm ,} the metric induced by $M$ on the $G$-orbits is
given by a bi-invariant pseudoRiemannian metric on $G$ that does
not depend on the $G$-orbit.
\end{prop}

\Proof Everywhere local freeness follows from topological
transitivity by the results in~\cite{Szaro}.

Observe that the condition for $G$-orbits to be nondegenerate is
an open condition, i.e.~there exist a $G$-invariant open subset
$U$ of $M$ so that the $G$-orbit of every point in $U$ is
nondegenerate.

On the other hand, given local freeness, it is well known that for
$T\OO$ the tangent bundle to the $G$-orbits, the following map is
a $G$-equivariant smooth trivialization of $T\OO$:
\begin{eqnarray*}
  \fhi : M\times \g &\rightarrow & T\OO \\
  (x,X) &\mapsto & X_x^*
\end{eqnarray*}
where $X^*$ is the vector field on $M$ whose one parameter group
of diffeomorphisms is $\exp(tX)$, and the $G$-action on
$M\times\g$ is given by $g(x,X) = (gx,{\rm Ad}(g)(X))$. Then, by
restricting the metric on $M$ to $T\OO$ and using the above
trivialization, we obtain the smooth map:
\begin{eqnarray*}
\psi : M &\rightarrow & \g^*\otimes\g^* \\
  x &\mapsto & B_x
\end{eqnarray*}
where $B_x(X,Y) = h_x(X_x^*,Y_x^*)$, for $h$ the metric on $M$.
This map is clearly $G$-equivariant. Hence, since the $G$-action
is tame on $\g^*\otimes\g^*$, such map is essentially constant on
the support of almost every ergodic component of $M$. Hence, if
$S$ is the support of one such ergodic component of $M$, then
there is an ${\rm Ad}(G)$-invariant bilinear form $B_S$ on $\g$ so
that, by the previous discussion, the metric on $T\OO|_S\cong
S\times\g$ induced by $M$ is almost everywhere given by $B_S$ on
each fiber. Also, the ${\rm Ad}(G)$-invariance of $B_S$ implies that
its kernel is an ideal of $\g$. If such kernel is $\g$, then
$T\OO|_S$ is lightlike which implies $\dim\g \leq m_0$. But this
contradicts the condition $n_0 = m_0$ since $n_0 < \dim\g$. Hence,
being $\g$ simple, it follows that $B_S$ is nondegenerate, and so
almost every $G$-orbit contained in $S$ is nondegenerate. Since
this holds for almost every ergodic component, it follows that
almost every $G$-orbit in $M$ is nondegenerate. In particular, the
set $U$ defined above is conull and so nonempty.

Moreover, the above shows that the image under $\psi$ of a conull,
and hence dense, subset of $M$ lies in the set of ${\rm Ad}(G)$-invariant
elements of $\g^*\otimes\g^*$. Since the latter set is closed, it
follows that $\psi(M)$ lies in it. In particular, on every $G$-orbit
the metric induced from that of $M$ is given by an
${\rm Ad}(G)$-invariant symmetric bilinear form on $\g$.

By topological transitivity, there is a $G$-orbit $\OO_0$ which is
dense and so it must intersect $U$. Since $U$ is $G$-invariant it
follows that $\OO_0$ is contained in $U$. Let $B_0$ be the
nondegenerate bilinear form on $\g$ so that under the map $\psi$
the metric of $M$ restricted to $\OO_0$ is given by $B_0$. Hence
$\psi(\OO_0) = B_0$ and so the density of $\OO_0$ together with
the continuity of $\psi$ imply that $\psi$ is the constant map
given by $B_0$. We conclude that all $G$-orbits are nondegenerate
as well as the last claim in the statement. \Endproof

The arguments in Proposition \ref{local-free} allows us to prove
the following result which is a generalization of Lemma
\ref{signatures}.

\begin{lemma}
\label{signatures2} Let $G$ be a connected noncompact simple Lie
group acting by isometries on a finite volume pseudoRiemannian
manifold $X$. Denote with $(n_1,n_2)$ and $(m_1,m_2)$ the
signatures of $G$ and $X${\rm ,} respectively{\rm ,} where $G$ carries a
bi-invariant pseudoRiemannian metric. If we denote $n_0 =
{\rm min}(n_1,n_2)$ and $m_0 = {\rm min} (m_1,m_2)${\rm ,} then $n_0 \leq m_0$.
\end{lemma}

\Proof With this setup we have local freeness on an open subset
$U$ of $X$ by the results in \cite{Zimmer-locally-free}. As in the
proof of Proposition \ref{local-free}, we consider the map:
\begin{eqnarray*}
U &\rightarrow & \g^*\otimes\g^* \\
  x &\mapsto & B_x
\end{eqnarray*}
which, from the arguments in such proof, is constant on the
ergodic components in $U$ for the $G$-action. On any such ergodic
component, the metric along the $G$-orbits comes from an
${\rm Ad}(G)$-invariant bilinear form $B_0$ on $\g$. As before, the
kernel of $B_0$ is an ideal. If the kernel is all of $\g$, then
$B_0 = 0$ and the $G$-orbits are lightlike which implies that $n_0
< \dim\g \leq m_0$. If the kernel is trivial, then $B_0$ is
nondegenerate and the $G$-orbits are nondegenerate submanifolds of
$X$. But this implies $n_0 \leq m_0$ as well, since $n_0$ does not
depend on the bi-invariant metric on $G$. \Endproof

In the rest of this work we will denote with $T\OO$ the tangent
bundle to the orbits. From Proposition \ref{local-free} it follows
that $TM = T\OO\oplus T\OO^\perp$, when the $G$-action is
topologically transitive and $n_0 = m_0$.

We will need the following result which provides large local
isotropy groups. Its proof relies heavily on the arguments in
\cite{GCT} (see also \cite{Gromov-rigid}). Similar results appear
in \cite{Zimmer-aut} and \cite{Feres}, but in such works
analyticity and compactness of the manifold acted upon is assumed.

\begin{prop}
\label{local-g} Let $G$ be a connected noncompact simple Lie group
and $X$ a smooth finite volume pseudoRiemannian manifold. Suppose
that $G$ acts smoothly on $X$ by isometries. Then there is a dense
subset $S$ of $X$ so that{\rm ,} for every $x\in S${\rm ,} there exist an open
neighborhood $U_x$ of $x$ and a Lie algebra $\g(x)$ of Killing
vector fields defined on $U_x$ satisfying\/{\rm :}\/
\begin{itemize}
    \item[{\rm (1)}] $Z_x = 0${\rm ,} for every $Z\in\g(x)${\rm ,}
    \item[{\rm (2)}] the local one-parameter subgroups of $\g(x)$
    preserve the $G$-orbits{\rm ,}
    \item[{\rm (3)}] $\g(x)$ and $\g$ are isomorphic Lie algebras{\rm ,} and
    \item[{\rm (4)}] for the isomorphism in {\rm (3),} the canonical vector space
    isomorphism\break $T_x Gx \cong \g$ is also an isomorphism of
    $\g$-modules.
\end{itemize}
\end{prop}

{\it Proof}. Without using analyticity, the arguments in Lemma 9.1 in
\cite{GCT} provide a conull set $S_0$ so that for every $x\in S_0$
one has a Lie algebra of infinitesimal Killing vector fields of
order $k$ that satisfy the above conclusions up to order $k$.
Moreover, this is achieved for every $k$ sufficiently large.
Further on, in Theorem 9.2 in \cite{GCT}, such infinitesimal
vector fields are extended to local ones by using analyticity.
This extension ultimately depends on Proposition~6.6 in
\cite{GCT}. The latter result is based on the arguments in Nomizu
\cite{Nomizu}.

In \cite{Nomizu} a notion of regular point for $X$ is introduced,
which satisfy two key properties. The set of regular points is an
open dense subset $U$ of $X$ and at regular points every
infinitesimal Killing field of large enough order can be extended
locally. The first property is found in \cite{Nomizu} and the
second one is proved in \cite{GCT}, both just using smoothness.
From these remarks we find that the set $S = U\cap S_0$ satisfies
the conclusions without the need to assume analyticity, as one
does in the statements of \cite{GCT}. Also, $S$ is obviously dense
since $U$ is open dense, $S_0$ is conull and the measure
considered (the pseudoRiemannian volume) is smooth. Finally, we
observe that even though the results in \cite{Nomizu} are stated
for Riemannian metrics, those that we use here extend with the
same proof to pseudoRiemannian manifolds. A remark of this sort
was already made in \cite{GCT}. \Endproof

We now prove integrability of the normal bundle to the orbits.

\begin{prop}
\label{integrable} Suppose $G$ acts topologically transitively on
$M$ preserving its pseudoRiemannian metric and satisfying $n_0 =
m_0$. Then $T\OO^\perp$ is integrable.
\end{prop}
\Proof Let $\omega : TM \rightarrow \g$ be the $\g$-valued 1-form
on M given by $TM = T\OO\oplus T\OO^\perp \rightarrow T\OO \cong
M\times \g \rightarrow \g$, where the two arrows are the natural
projections. Define the curvature of $\omega$ by the 2-form
$\Omega = d\omega|_{T\OO^\perp\wedge T\OO^\perp}$. As remarked in
\cite{Gromov-rigid} (see also \cite{Paco}) it is easy to prove
that $T\OO^\perp$ is integrable if and only if $\Omega = 0$.

Choose $S$ and $\g(x)$ as in Proposition \ref{local-g}. Hence, the
local one-parameter subgroups of $\g(x)$ preserve $T\OO^\perp_x$
for every $x \in S$. From this, and the isomorphism $\g(x)\cong
\g$ described in the proof of Proposition \ref{local-g}, it is
easy to show that the linear map $\Omega_x : T\OO^\perp_x\wedge
T\OO^\perp_x \rightarrow \g$ is a $\g$-module homomorphism, for
every $x \in S$. This fact is contained in the proof of
Proposition 3.9 in \cite{Paco}.

On the other hand, Proposition \ref{local-free} and the condition
$n_0 = m_0$ imply that $T\OO^\perp$ is either Riemannian or
antiRiemannian. Since the elements of $\g(x)$ are Killing fields,
it follows that $\g(x)$ can be linearly represented on
$T\OO^\perp_x\wedge T\OO^\perp_x$ so that the elements of $\g(x)$
define derivations of a definite inner product. This provides a
homomorphism of $\g(x)$ into the Lie algebra of an orthogonal
group. Since $\g(x)$ is simple noncompact, such homomorphism is
trivial and it follows that the $\g(x)$-module $T\OO^\perp_x\wedge
T\OO^\perp_x$ is trivial. Hence no $\g(x)$-irreducible subspace of
$T\OO^\perp_x\wedge T\OO^\perp_x$ can be isomorphic to $\g$ as a
$\g$-module. By Schur's Lemma we conclude that $\Omega_x = 0$, for
every $x\in S$. Hence, $\Omega = 0$ on all of $M$ and $T\OO^\perp$
is integrable. \Endproof

The following result is fundamental to obtain an isometric
splitting.

\begin{lemma}
\label{complete-N} With the conditions of Proposition
{\rm \ref{integrable},} the leaves of the foliation defined by
$T\OO^\perp$ are complete for the metric induced by $M$.
\end{lemma}
\Proof As observed in the proof of Proposition \ref{integrable},
the bundle $T\OO^\perp$ is either Riemannian or antiRiemannian.
Hence, the foliation by $G$-orbits on $M$ carries a Riemannian or
antiRiemannian structure obtained from $T\OO^\perp$. By the basic
properties of Riemannian foliations, the compactness of $M$
implies that geodesic completeness is satisfied for geodesics
orthogonal to the $G$-orbits (see \cite{Molino}). This clearly
implies the completeness for leaves of the foliation given by
$T\OO^\perp$. \Endproof

The next proposition provides a first description of the
properties of $M$. It is similar in spirit to the main results in
\cite{Paco}.

\begin{prop}
\label{split1} Suppose $G$ acts topologically transitively on $M$
preserving its pseudoRiemannian metric and satisfying $n_0 = m_0$.
Choose a leaf $N$ of the foliation defined by $T\OO^\perp$. Fix on
$G$ the bi-invariant pseudoRiemannian metric that induces on the
$G$-orbits the metric inherited by $M$ and consider $N$ as a
pseudoRiemannian manifold with the metric inherited by $M$ as
well. Then the map $G\times N \rightarrow M${\rm ,} obtained by
restricting the $G$-action\break to $N${\rm ,} is a $G$-equivariant isometric
covering map. In particular{\rm ,} this induces a $G$-equivariant
isometric covering map $G\times \widetilde{N} \rightarrow M${\rm ,}
where $\widetilde{N}$ is the universal covering space of $N$.
\end{prop}

\vglue-3pt {\it Proof}.  By our choices of metrics, the $G$-invariance of the metric
on $M$ and the previous results, it is easy to conclude that the
map $G\times N \rightarrow M$ as above is a local isometry. On the
other hand, the Levi-Civita connection on $G$ is bi-invariant and,
by the problems in Chapter II of \cite{Helgason}, its geodesics
are translates of one-parameter subgroups. In particular, $G$ is
complete. Hence, by Lemma \ref{complete-N}, $G\times N$ is a
complete pseudoRiemannian manifold. Then, Corollary 29 in page 202
in \cite{ONeill} implies that the restricted action map $G\times N
\rightarrow M$ is an isometric covering map. The rest of the
claims follow easily from this fact. \phantom{overthere}\Endproof\vskip4pt

As an immediate consequence we obtain the following result.  The proof uses Proposition 4.5.  We note that
Section~4 is actually independent from this section and the rest of this work.

\begin{cor}
\label{lattice1} Let $G$, $M$ and $N$ be as in the hypotheses of
Proposition~{\rm \ref{split1}.} Then there is a discrete subgroup
$\Gamma_0$ of $\Iso(G\times \widetilde{N})$ of deck transformations for $G\times \widetilde{N}\to M$
such that
$ (G\times \widetilde{N})/\Gamma_0 \to M$ is a $G$-equivariant finite covering.
\end{cor}

Our next goal is to prove that, by passing to a finite covering,
$\Gamma_0$ can be replaced by a discrete subgroup of a group that
contains $G$ as a subgroup as well. In order to do that, we will
study the isometry group of $G$ with some bi-invariant
pseudoRiemannian metric.

\section{Geometry of bi-invariant metrics on $G$}

\label{section-G} Let $G$ be a connected noncompact simple Lie group
$G$ as before. We will investigate some useful properties about the
geometry of a bi-invariant pseudoRiemannian metric on $G$. Note that
any such metric is analytic and by \cite{Paco-thesis} can be
described in terms of the Killing form (see also \cite{Bader-Nevo});
however, we will not use such fact. In this section, we fix an
arbitrary bi-invariant pseudoRiemannian metric on $G$ and denote
with $\Iso(G)$ the corresponding group of isometries. Also we denote
$L(G)R(G) = \{ L_g\circ R_h | g,h \in G \}$, the group generated by
left and right translations, which is clearly a connected subgroup
of $\Iso(G)$.

We will use some basic properties of pseudoRiemannian symmetric
spaces, which are known to be a natural generalization of
Riemannian symmetric spaces. For the definitions and basic
properties of the objects involved we will refer to
\cite{Cahen-Parker}. Moreover, we will use in our proofs some of
the results found in this reference.

From \cite{Cahen-Parker}, we recall that, in a pseudoRiemannian
symmetric space $X$, a transvection is an isometry of the form
$s_x\circ s_y$, where $s_x$ is the involutive isometry that has
$x$ as an isolated fixed point. The group $T$ generated by
transvections is called the transvection group of $X$. This group
is (clearly) invariant under the conjugation by $s_o$, any fixed
involutive isometry. With this setup, the pseudoRiemannian
symmetric triple associated to $X$ is given by
$(Lie(T),\sigma,B)$, where $B$ is a suitable bilinear form on
$Lie(T)$ and $\sigma$ is the differential at $e\in T$ of the
conjugation by some involution $s_o$. We refer to
\cite{Cahen-Parker} for a more precise description of this object.
Here, we need to show the following features of the geometry of
$G$ associated with these notions.

\begin{prop}
\label{prst-G} $G$ is a pseudoRiemannian symmetric space whose
associated pseudoRiemannian symmetric triple can be chosen to be
of the form $(\g\times\g,\sigma,B)${\rm ,} where $\sigma(X,Y) = (Y,X)$.
\end{prop}
\Proof Since the differential of the inversion map $g \mapsto
g^{-1}$ at any point can be written as the composition of the
differentials of a left and a right translations (see
\cite{Helgason}), it follows that the bi-invariant metric on $G$
is also invariant under the inversion map. Hence, for every $x\in
G$, the map $s_x$ defined by $s_x(g) = xg^{-1}x$ is an isometry of
$G$ and it is easily seen to be involutive with $x$ as an isolated
fixed point. Hence, $G$ is pseudoRiemannian symmetric.

Let $T$ be the transvection group of $G$. One can easily check
that $s_x\circ s_y = L_{xy^{-1}}\circ R_{y^{-1}x}$, and so $T$ is
a subgroup of $L(G)R(G)$. On the other hand, since $G$ is simple
and connected we have $[G,G] = G$. Hence, for every $z\in G$,
there exist $x,y\in G$ such that $z = [x,y]$. From this it is easy
to prove that $L_z = s_e\circ s_{yx}\circ s_x\circ s_{y^{-1}}\in
T$. This with a similar formula for right translations show that
$T = L(G)R(G)$. Furthermore, if we define a $G\times G$-action on
$G$ by $(g,h)x = gxh^{-1}$, then the map $(g,h) \mapsto L_g\circ
R_{h^{-1}}$ defines a local isomorphism $G\times G \rightarrow
L(G)R(G) = T$, which implies $Lie(T) \cong \g \times \g$ as Lie
algebras. Finally, a straightforward computation proves that using
conjugation by $s_e$, the map $\sigma$ on $\g \times \g$ has the
required expression. \Endproof

As a consequence we obtain the following result. We recall that a
connected pseudoRiemannian manifold is called weakly irreducible if
the tangent space at some (and hence at every) point has no
nonsingular proper subspaces invariant by the holonomy group at the
point.

\begin{prop}
\label{G-w-irred} For any bi-invariant pseudoRiemannian metric on
$G${\rm ,} the universal covering space $\widetilde{G}$ is weakly
irreducible.
\end{prop}
\Proof Consider the representation $\rho$ of the Lie algebra $\f =
\{ (X,X) | X\break \in \g \}$ (isomorphic to $\g$) in the space $\p =
\{(Y, -Y) | Y \in \g \}$ given by the expression:
$$
        \rho(X,X)(Y,-Y) = ([X,Y],-[X,Y]).
$$
This clearly turns $\p$ into an $\f$-module isomorphic to the
$\g$-module given by the adjoint representation of $\g$. Since
$\g$ is simple, $\p$ is then an irreducible $\f$-module. Then the
conclusion follows from the description of the pseudo\-Riemannian
symmetric triple associated to $G$ in our Proposition \ref{prst-G}
and from Proposition~4.4 in page 18 in \cite{Cahen-Parker}.
\Endproof
\vskip5pt

With the previous result at hand we obtain the next statement.

\begin{prop}
\label{iso-GxN} Let $N$ be a connected complete Riemannian \/{\rm (}\/or
antiRiemannian\/{\rm )}\/ manifold. Then any isometry of the
pseudoRiemannian product $G\times N$ preserves the factors{\rm ,} in
other words{\rm ,} $\Iso(G\times N) = \Iso(G)\times\Iso(N)$.
\end{prop}
\Proof Let $\widetilde{G}\times \widetilde{N}$ be the universal
covering of $G\times N$. Let $\widetilde{N} = N_0\times\dots
\times N_k$ be the de Rham decomposition of $\widetilde{N}$ as
Riemannian (or antiRiemannian) manifold. By the de Rham-Wu
decomposition theorem for pseudoRiemannian manifolds (see
\cite{Wu} and \cite{Cahen-Parker}) and by Proposition
\ref{G-w-irred}, it follows that $\widetilde{G}\times
\widetilde{N}$ has a de Rham decomposition and it is given by
$\widetilde{G}\times N_0\times\dots N_k$. Furthermore, it is known
that such decomposition is unique up to order. In particular,
every isometry of $\widetilde{G}\times N_0\times\dots N_k$
permutes the factors, but since each $N_i$ is Riemannian (or
antiRiemannian) and $\widetilde{G}$ is not, then every isometry of
$\widetilde{G}\times \widetilde{N}$ preserves these two factors.

Now let $f\in\Iso(G\times N)$ and lift it to an isometry
$\widetilde{f}$ of $\widetilde{G}\times \widetilde{N}$. By the
previous arguments, $\widetilde{f}$ preserves the product, i.e. if
we write $\widetilde{f} = (\widetilde{f}_1,\widetilde{f}_2)$, then
$\widetilde{f}_1$ does not depend on $\widetilde{N}$ and
$\widetilde{f}_2$ does not depend on $\widetilde{G}$. From this,
it is easy to see that there exist isometries $f_1\in\Iso(G)$ and
$f_2\in\Iso(N)$ such that $f = (f_1,f_2)$, thus showing the
result. \Endproof
\vskip5pt

We now proceed to obtain a fairly precise description of
$\Iso(G)$. First we prove the following result.

\begin{lemma}
\label{IsoGe-closed} Denote by $\Iso(G)_e$ the isotropy subgroup
at $e\in G$. Then the homomorphism\/{\rm :}\/
\begin{eqnarray*}
\fhi : \Iso(G)_e & \rightarrow & \Gl(\g) \\
    h & \mapsto & dh_e
\end{eqnarray*}
is an isomorphism onto a closed subgroup of $\Gl(\g)$.
\end{lemma}
\Proof By the arguments from Lemma 11.2 in page 62 in
\cite{Helgason} the map is injective. Now let $L^{(1)}(G)$ be the
linear frame bundle of $G$ endowed with the parallelism given by
the Levi-Civita connection on $G$. Consider the standard fiber of
$L^{(1)}(G)$ given by $\Gl(\g)$. Then the proof of Theorem 3.2 in
page 15 in \cite{Kobayashi} shows that the map:
\begin{eqnarray*}
 \psi : \Iso(G) &\rightarrow &  L^{(1)}(G) \\
  h &\mapsto & dh_e
\end{eqnarray*}
realizes $\Iso(G)$ as a closed submanifold that defines the
topology for which $\Iso(G)$ turns out to be a transformation Lie
group. Then the result follows by observing that the image of
$\fhi$ is just the intersection of the image of $\psi$ with the
fiber of $L^{(1)}(G)$ at $e$, which is $\Gl(\g)$. \Endproof\vskip5pt

The next result turns out to be fundamental in our proof of
Theorem A.

\begin{prop}
\label{IsoG-component} $\Iso(G)$ has finitely many components and
$\Iso(G)_0 = L(G)R(G)$. Also{\rm ,} $(\Iso(G)_e)_0$ is isomorphic to
${\rm Ad}(G)$ with respect to the map $\fhi $ from Lemma
{\rm \ref{IsoGe-closed}.}
\end{prop} 

\Proof By Problem A.6.ii in page 148 in \cite{Helgason} and since
$\Iso(G)_e$ preserves the (unique) bi-invariant connection we have
that $\fhi(\Iso(G)_e) \subset \Aut(\g)$, where $\Aut(\g)$ is the
group of Lie algebra automorphisms of $\g$. Also, we clearly have
${\rm Ad}(G) \subset \fhi(\Iso(G)_e)$, because $L_g\circ R_{g^{-1}}$ is
an isometry. Since $\g$ is simple we know that ${\rm Ad}(G) =
\Aut(\g)_0$, from which we conclude $\fhi(\Iso(G)_e)_0 = {\rm Ad}(G)$,
and the last claim follows.

Since $\Aut(\g)$ is algebraic, it has finitely many components.
Hence, the previous inclusions imply that $\Iso(G)_e$ has finitely
many components and, by Lemma \ref{finite-components}, the group
$\Iso(G)$ satisfies such property as well.

On the other hand, $G$ being homogeneous under $\Iso(G)$ it is
also homogeneous under $\Iso(G)_0$ (see \cite{Helgason}). Hence,
we observe that $G$ has the following two expressions as a
homogeneous pseudoRiemannian manifold:
\begin{eqnarray*}
  G &=& \Iso(G)_0/(\Iso(G)_0\cap\Iso(G)_e) \\
    &=& L(G)R(G)/G_\Delta,
\end{eqnarray*}
where $G_\Delta= \{ L_g\circ R_{g^{-1}} | g \in G \}$. This
together with the following easy to prove identities:
$$
\dim(\Iso(G)_0\cap\Iso(G)_e) = \dim((\Iso(G)_e)_0) = \dim({\rm Ad}(G))
= \dim(G_\Delta),
$$
show that the inclusion $L(G)R(G) \subset \Iso(G)_0$ is a
homomorphic inclusion of connected Lie groups of the same
dimension. This implies $\Iso(G)_0 = L(G)R(G)$. \Endproof\vskip5pt

By Propositions \ref{iso-GxN}, \ref{IsoG-component} and Corollary
\ref{isoX} we obtain the following.

\begin{prop}
\label{iso-GxN-component} Let $N$ be a connected homogeneous
Riemannian manifold. Then $\Iso(G\times N)$ has finitely many
connected components and\break $\Iso(G\times N)_0 =
\Iso(G)_0\times\Iso(N)_0 = L(G)R(G)\times\Iso(N)_0$.
\end{prop}

\section{Proof of Theorem A}
\label{section-thmA} Throughout this section, we assume that $G$ and
$M$ satisfy the hypotheses of Theorem A. Let $N$ be a leaf of the
foliation defined by $T\OO^\perp$, as described in
Section~\ref{section-split}, and $\widetilde{N}$ its universal
covering space. We will prove that $\widetilde{N}$ is homogeneous.
First observe that, by hypothesis, there is a dense $G$-orbit. By
Proposition \ref{split1}, we have $GN = M$, and so we can assume
that for some $x_0\in N$, the $G$-orbit of $x_0$ is dense in $M$.

Let $y$ be any given point in $N$ and $V$ an open neighborhood of
$y$ in $N$ with the leaf topology. Since the complementary
distributions $T\OO$ and $T\OO^\perp$ (as defined in \S
\ref{section-split}) are integrable, we can find a connected
neighborhood $U$ of $y$ in $M$ and a diffeomorphism $f = (f_1,f_2)
: U \rightarrow \R^k\times\R^l$ that maps the pair of foliations
for the distributions $T\OO|_U$ and $T\OO^\perp|_U$
diffeomorphically onto the foliations of $\R^{k+l}$ defined by the
factors $\R^k$ and $\R^l$ (see Example 2.10 in page~12 in
\cite{Kobayashi}). By shrinking either $U$ or $V$, if necessary,
we can assume that $V$ is a plaque for the foliated chart that $f$
defines for the foliation given by $T\OO^\perp$. Then $V$ is a
transversal for the foliation by $G$-orbits that intersects every
plaque (for the foliation by $G$-orbits) in $U$. Let $x$ be a
point in the $G$-orbit of $x_0$ lying in $U$, and $P$ the plaque
of the foliation by $G$-orbits that contains $x$. If we choose
$z\in P\cap V$, since the plaque $P$ is contained in the $G$-orbit
of $x_0$ (because $x\in P$), then there exist $g\in G$ such that
$z = g x_0$. But the $G$-action preserves $T\OO^\perp$, and so the
restriction of $g$ to $N$ defines an isometry of $N$ that maps
$x_0$ into $z$.

The previous arguments show that the isometry group of $N$ has a
dense orbit in $N$. By lifting isometries of $N$ to $\widetilde{N}$
and using deck transformations, it is easy to prove that the
isometry group of $\widetilde{N}$ has a dense orbit as well. Since
$\widetilde{N}$ is either Riemannian or antiRiemannian (see the
proof of Proposition~\ref{integrable}), we conclude from Lemma
\ref{complete-N} and Theorem \ref{singer-thm} that $\widetilde{N}$
is a homogeneous pseudoRiemannian manifold. Let $H$ be the identity
component of the group of isometries of $\widetilde{N}$. Then there
is a compact subgroup $K$ of $H$ such that
$\widetilde{N}=K\backslash H$.

Let $\Gamma_0$ be a discrete subgroup $\Iso(G\times \widetilde{N})$
as given by Corollary \ref{lattice1}. By Proposition
\ref{iso-GxN-component}, the group $\Iso(G\times \widetilde{N})$ has
finitely many components and its identity component is
$L(G)R(G)\times\Iso(N)_0 = L(G)R(G)\times H$. If we define $\Gamma =
\Gamma_0\cap (L(G)R(G)\times H)$, then $\Gamma$ is a normal finite
index subgroup of $\Gamma_0$, which provides a normal finite
covering $\widehat{M}= (G\times \widetilde{N})/\Gamma = (G\times
K\backslash H)/\Gamma$ of~$M$. Now, let $\gamma\in \Gamma$ be given.
Then, by the lifting property, we can find a diffeomorphism
$\widetilde{\gamma}$ such that the following diagram \pagebreak commutes.
$$
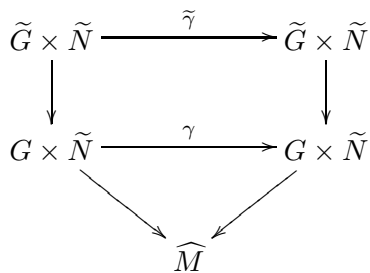
\xymatrix{ \widetilde{G}\times \widetilde{N}
\ar[rr]^{\widetilde{\gamma}} \ar[d]&    &
                    \widetilde{G}\times\widetilde{N} \ar[d] \\
G\times \widetilde{N} \ar[rr]^\gamma \ar[dr] &    &
                    G\times \widetilde{N} \ar[dl]  \\
&   \widehat{M} &
 }
$$
where the vertical arrows are given by the universal covering map
of $G\times \widetilde{N}$ and the arrows into $\widehat{M}$ are
given by the covering map $G\times \widetilde{N} \rightarrow
\widehat{M}$. In particular, $\widetilde{\gamma}$ is a covering
transformation for the universal covering map $\widetilde{G}\times
\widetilde{N} \rightarrow \widehat{M}$, and so it is also a
covering transformation for the universal covering map
$\widetilde{G}\times \widetilde{N} \rightarrow M$. By Proposition
\ref{split1}, the $G$-action on $M$ lifts to a $G$-action on
$G\times \widetilde{N}$ and to a $\widetilde{G}$-action on
$\widetilde{G}\times \widetilde{N}$ both by the expression
$g(g_1,x) = (gg_1,x)$. On the other hand, the latter lifted action
commutes with the action of $\pi_1(M)$ and in particular it
commutes with $\widetilde{\gamma}$. Then we observe that, in the
square of the commutative diagram above, the top horizontal arrow
and the vertical arrows are $\widetilde{G}$-equivariant and so the
map $\gamma$ is $G$-equivariant. On the other hand, we can write
$\gamma = (L_{g_1}\circ R_{g_2},h)$, since it lies in
$L(G)R(G)\times H$. Then, the $G$-equivariance of $\gamma$ and the
above expression for the $G$-action on $G\times \widetilde{N}$
yields:
$$
    (g_1gg_2,xh) = (g(e,x))(L_{g_1}\circ R_{g_2},h) =
    g((e,x)(L_{g_1}\circ R_{g_2},h)) = (gg_1g_2,xh)
$$
for every $g\in G$ and $x\in \widetilde{N}$, which implies $g_1
\in Z(G)$. Hence, $L_{g_1} = R_{g_1}$ and then $\gamma \in
R(G)\times H$.

The previous arguments show that $\Gamma$ is a discrete subgroup
of $R(G)\times H \cong G\times H$ acting on the right on $G\times
\widetilde{N} = G\times K\backslash H$ and commuting with the
(natural) $G$-action on the latter space. In particular, $\Gamma$
is a discrete cocompact subgroup of $G\times H$. Since the
$\Gamma$-action is on the right and the $G$-action is on the left,
they both commute and the $G$-action lifts to $\widehat{M}$.

If we now take $L = G\times H$, $\widehat{M}$, $\Gamma$ and $K$ as
above, then the conditions required in (1), (2) and (3) are
satisfied, as well as the fact that $\widehat{M}$ is\break
$G$-equivariantly diffeomorphic to $K\backslash L/\Gamma$. Hence,
to complete the proof of Theorem A, it remains to show that the
$G$-actions on $M$ and on $L/\Gamma$ are ergodic and engaging. We
will prove this by showing that we can replace $H$, $\Gamma$ and
$\widehat{M}$ by suitable choices so that the properties proved so
far still hold and the required ergodicity and engagement
conditions are now satisfied. To achieve this, we will study the
properties of the Riemannian (or antiRiemannian) foliation by
$G$-orbits on $M$.

We will use the following result about Riemannian and Lie
foliations.

\begin{lemma}
\label{lie-riem-foliations} Let $X$ be a compact manifold with a
smooth foliation $\FF$ carrying either a transverse Lie structure
or a transverse Riemannian structure. Then the foliation $\FF$ has
a dense leaf if and only if it is ergodic. Furthermore{\rm ,} if $\FF$
has a dense leaf{\rm ,} then the induced foliation on any connected
finite covering of $X$ has a dense leaf as well and so it is
ergodic.
\end{lemma}
\Proof In both cases, the foliation carries a transverse holonomy
invariant finite volume and so ergodicity is easily seen to imply
the existence of a dense leaf.

For the converse, let us first assume that $\FF$ carries a
transverse Lie structure. Let $F$ be a connected Lie group that
models the transverse Lie structure, so that there exist a
developing map and a holonomy representation given by:
\begin{eqnarray*}
  D : \widehat{X} &\rightarrow& F \\
  \rho : \pi_1(X) &\rightarrow& F
\end{eqnarray*}
We refer to \cite{Molino}, \cite{Tondeur} and \cite{Zimmer-Lie}
for the definition of such objects.

Then it is easily seen that the existence of a dense leaf is
equivalent to the density of $\rho(\pi_1(X))$ in $F$. By Lemma
2.2.13 in \cite{Zimmer-book} such density is equivalent to the
ergodicity of the $\rho(\pi_1(X))$-action on $F$ by left
translations. Finally, the latter ergodicity is equivalent to the
ergodicity of the foliation $\FF$ (see for example \cite{FHM} or
\cite{Zimmer-Lie}).

We still consider $X$ endowed with a transverse Lie structure and
now we suppose that $X$ has a dense leaf. Let $Y \rightarrow X$ be a
connected finite covering map. If we choose a dense leaf $\OO$ in
$X$, then, for the induced foliation on $Y$, the inverse image of
$\OO$ in $Y$ is a finite union of leaves $\OO_1, \dots, \OO_r$. We
clearly have that $\OO_1\cup \dots\cup \OO_r$ is dense in $Y$, i.e.
we can write $Y = \overline{\OO}_1\cup \dots\cup \overline{\OO}_r$.
On the other hand, since $Y$ is compact, the transverse Lie
structure on $Y$ is complete and by the structure theorems for
transversely parallelizable foliations (see \cite{Molino}) it
follows that the closures of leaves define a partition of $Y$. Since
$Y$ is connected, from the previous expression of $Y$ it follows
that $Y = \overline{\OO}_i$, for some $i$, and so $Y$ has a dense
leaf. Moreover, from what we have proved so far it follows that the
foliation in $Y$ is ergodic.

We now consider the case where $X$ carries a transverse Riemannian
structure. Then, by the results in \cite{Molino}, there is a
foliated fiber bundle $\pi : P \rightarrow X$ where $P$ is compact
and carries a transverse Lie structure. Moreover, if $X$ has a dense
leaf it is easy to see from the results in \cite{Molino} (see page
153 therein) that we can choose a foliated reduction to assume that
$P$ has a dense leaf. In particular, from the above, the foliation
on $P$ is ergodic. Hence, if the foliation $\FF$ in $X$ is not
ergodic, then there exist a leaf saturated measurable subset $A$ of
$X$ which is neither null nor conull. From this, it is easily seen
that $\pi^{-1}(A)$ is a measurable leaf saturated subset of $P$
which is neither null nor conull. This contradiction proves that
$\FF$ is ergodic.

Finally, the last claim for transverse Riemannian structures is
proved with the same argument used for Lie foliations using the
corresponding properties for Riemannian foliations. \hfill\qed\vskip5pt

As an immediate consequence of Lemma \ref{lie-riem-foliations},
the $G$-action on $M$ is ergodic and engaging.

We now state the following result. Its proof follows easily from
the definitions, and the (easy) fact that there is an isomorphism
\smallbreak
\centerline{$\pi_1((G\times H)/\Gamma)/\pi_1(G\times H) \cong \Gamma$}
\smallbreak\nobreak\noindent (see
\cite{Hu} or \cite{KN}).

\begin{lemma}
\label{developing} Let $G$, $H$ and $\Gamma$ be given as above.
Let $D : G\times H \rightarrow H$ be the natural projection.
Define $\rho : \pi_1((G\times H)/\Gamma) \rightarrow H$ as the
composition of the maps\/{\rm :}\/
$$
\pi_1((G\times H)/\Gamma) \rightarrow \pi_1((G\times
H)/\Gamma)/\pi_1(G\times H) \cong \Gamma \rightarrow H
$$
where the last arrow is the restriction of $D$ to $\Gamma$. Then
$D$ is a developing map with holonomy representation $\rho$ for
the Lie foliation on $(G\times H)/\Gamma$ given by the local
factor $H$.
\end{lemma}

We observe that $G$ preserves the (finite) pseudoRiemannian volume
on $M$ and so it preserves a finite smooth measure on $K\backslash
L/\Gamma$. Since $K$ is compact, it is easy to check that there is a
$G$-invariant finite smooth measure on $L/\Gamma$. Moreover, from
the proof of Lemma \ref{lie-riem-foliations}, it is easy to see that
the $G$-action on $L/\Gamma$ is ergodic (with respect to the latter
$G$-invariant measure) if and only if ${\rm Im}(\rho)$ is dense in $H$. If
such density holds with our current choices of groups, then by Lemma
\ref{lie-riem-foliations} the $G$-action on $L/\Gamma$ is ergodic
and engaging.

On the other hand, if ${\rm Im}(\rho)$ is not dense in $H$, we will show
that we can replace $H$ by a smaller connected closed subgroup
$H_1$ and $\Gamma$ by a normal finite index subgroup $\Gamma_1$
for which the conclusions of Theorem $A$ proved up to this point
still hold and such that the corresponding holonomy representation
has dense image in $H_1$.

\begin{lemma}
\label{H-reduction} Let $P_T$ be the orthogonal transverse frame
bundle of the\break Riemannian \/{\rm (}\/or antiRiemannian\/{\rm )}\/ foliation on
$\widehat{M}$. Then the $G$-action lifts to $P_T$ and there is a
$G$-equivariant map $\fhi : (G\times H)/\Gamma \rightarrow P_T$
for which the diagram\/{\rm :}\/
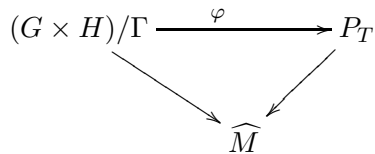
$$
\xymatrix{ (G\times H)/\Gamma \ar[rr]^\fhi \ar[dr] &    &
                    P_T \ar[dl]  \\
&   \widehat{M} &
 }
$$
is commutative and realizes $(G\times H)/\Gamma$ as an embedded
$K$-reduction of $P_T$. Here we recall that $\widetilde{N} =
K\backslash H$.
\end{lemma}
\Proof Since $\widehat{M} = (G\times \widetilde{N})/\Gamma$, it is
a simple matter to prove that the transverse frame bundle over
$\widehat{M}$ associated to the foliation by $G$-orbits is given
by $(G\times L^{(1)}(\widetilde{N}))/\Gamma$, where
$L^{(1)}(\widetilde{N})$ is the linear frame bundle of
$\widetilde{N}$. Note that in our notation the group of
diffeomorphisms on both $N$ and $\widetilde{N}$ has been chosen to
act on the right and so $L^{(1)}(\widetilde{N})$, and its
reductions, have left actions for their structure groups.

From the previous description, the $G$-action on $\widehat{M}$
lifts to the transverse frame bundle of the foliation on
$\widehat{M}$ by $G$-orbits. Hence, since $G$ preserves the
Riemannian (or antiRiemannian) structure on the foliation, then
the $G$-action preserves $P_T$. The latter action is thus the lift
of the $G$-action on $\widehat{M}$.

Let $o=Ke \in K\backslash H$ and consider $\Gl(T_o\widetilde{N})$
as the structure group of $L^{(1)}(\widetilde{N})$. Define the
map:
\begin{eqnarray*}
\widetilde{\fhi} : G\times H &\rightarrow& G\times L^{(1)}(\widetilde{N})  \\
  (g,h) &\mapsto& (g,dh_o)
\end{eqnarray*}
Such map is clearly equivariant with respect to both $G$ and
$\Gamma$. On the other hand, the restriction of $\widetilde{\fhi}$
to $\{ e\}\times K$ clearly realizes $K$ as a closed subgroup of
$O(T_o\widetilde{N})$. Also, since $H$ acts by isometries on
$\widetilde{N}$, it induces a $G$-equivariant map $\fhi : (G\times
H)/\Gamma \rightarrow P_T$.

On the other hand, the natural projection $(G\times H)/\Gamma
\rightarrow (G\times \widetilde{N})/\Gamma = \widehat{M}$ is
clearly surjective, because $\widetilde{N} = K\backslash H$, and
its fibers are $K$-orbits diffeomorphic to $K$ since $H$ acts
effectively on $\widetilde{N}$. It follows that, with respect to
$\fhi$, the manifold $(G\times H)/\Gamma$ is a $K$-reduction of
$P_T$. \Endproof

Now let $H$, $\Gamma$ and $\rho$ be as before, where the latter is
the holonomy map from Lemma \ref{developing}. Recall that we are
assuming that $\overline{{\rm Im}(\rho)}$ is a proper subgroup of $H$.
First we observe that by the definition of $\rho$, the group
$\Gamma$ is a subgroup of $G\times \overline{{\rm Im}(\rho)}$. Hence, we
clearly have that $(G\times \overline{{\rm Im}(\rho)})/\Gamma$ is a
closed $G$-invariant subset of $(G\times H)/\Gamma$. In
particular, $(G\times \overline{{\rm Im}(\rho)})/\Gamma$ is a union of
closures of $G$-orbits. By the structure of the closures of leaves
of a Riemannian foliation (see Theorem 5.1 in page 155 in
\cite{Molino}), we know that the closure of any $G$-orbit in $P_T$
is mapped onto the closure of some $G$-orbit in $\widehat{M}$.
Since we proved that the $G$-action on $M$ is engaging, we also
know that in $\widehat{M}$ every $G$-orbit is dense, in other
words, $\widehat{M}$ is a single leaf closure. Hence, any
$G$-invariant closed subset in $(G\times H)/\Gamma$ maps onto
$\widehat{M}$ under the projection $(G\times H)/\Gamma \rightarrow
(G\times K\backslash H)/\Gamma = \widehat{M}$. Since $(G\times
\overline{{\rm Im}(\rho)})/\Gamma$ is $G$-invariant and closed, it is
mapped onto $\widehat{M}$ under the latter natural projection.

\begin{lemma}
\label{H_1-homogeneity} With the above setup{\rm ,}
$\overline{{\rm Im}(\rho)}$ is a closed subgroup of
$\Iso(\widetilde{N})$ that acts transitively on $\widetilde{N}$.
In particular{\rm ,} if $H_1$ is the identity component of
$\overline{{\rm Im}(\rho)}${\rm ,} then $H_1$ acts transitively on
$\widetilde{N}$.
\end{lemma}

\Proof Let $o = Ke \in K\backslash H = \widetilde{N}$. Choose an
arbitrary $y\in \widetilde{N}$. Then by the previous discussion,
there exist $(g,h) \in G\times \overline{{\rm Im}(\rho)}$ and
$(\gamma_1,\gamma_2) \in \Gamma \subset G\times
\overline{{\rm Im}(\rho)}$ such that:
$$
(e,y) = (g,Kh)(\gamma_1,\gamma_2) = (g\gamma_1, Kh\gamma_2) =
(g\gamma_1, oh\gamma_2)
$$
and so $y = oh\gamma_2$, where $h\gamma_2 \in
\overline{{\rm Im}(\rho)}$. Hence, $\overline{{\rm Im}(\rho)}$ acts
transitively on $\widetilde{N}$ and, by connectedness of
$\widetilde{N}$, so does $H_1$. \Endproof

By Lemma \ref{H_1-homogeneity} and Corollary \ref{isoX}, it
follows that $H_1$ is a normal finite index subgroup of
$\overline{{\rm Im}(\rho)}$. Hence, the group $\Gamma_1 = (G\times
H_1)\cap\Gamma$ is a normal finite index subgroup of $\Gamma$. Let
$K_1$ be the stabilizer in $H_1$ of the point $o = Ke \in
K\backslash H = \widetilde{N}$. If we define the manifold
$\widehat{M}_1 = (G\times K_1 \backslash H_1)/\Gamma_1$, then,
from the previous constructions, the natural projection map:
$$
(G\times K_1 \backslash H_1)/\Gamma_1 \rightarrow (G\times K
\backslash H)/\Gamma = (G\times \widetilde{N})/\Gamma =
\widehat{M}
$$
defines a finite normal covering of $\widehat{M}$ and so
$\widehat{M}_1$ is a finite covering of $M$.

We claim that if we choose $\widehat{M}_1$ for (1), $L = G\times
H_1$ for (2) and $\Gamma_1 \subset L$ and $K_1 \subset C_L(G)$ for
(3), then the conclusions of Theorem A are satisfied. In fact, the
only point that remains to be proved is that the $G$-action on
$L/\Gamma_1$ is ergodic and engaging.

To show this last requirement, we first observe that Lemma
\ref{developing} still holds with our new choices. In other words,
the developing map $D_1 : G\times H_1 \rightarrow H_1$ and the
holonomy representation $\rho_1 : \pi_1((G\times H_1)/\Gamma_1)
\rightarrow H_1$ of the transverse Lie structure for the foliation
by $G$-orbits on $(G\times H_1)/\Gamma_1$ are the restrictions to
$G\times H_1$ and $\pi_1((G\times H_1)/\Gamma_1)$, respectively,
of the maps $D$ and $\rho$ as they are defined in Lemma
\ref{developing}.

On the other hand, from the proof of Lemma
\ref{lie-riem-foliations}, to conclude ergodicity and engagement
of the $G$-action on $L/\Gamma_1$ it is enough to show that
$\overline{{\rm Im}(\rho_1)} = H_1$. From the expressions of $\rho$ and
$\rho_1$ in Lemma \ref{developing}, if we denote with $pr :
G\times H \rightarrow H$ the natural projection, then ${\rm Im}(\rho) =
{\rm pr}(\Gamma)$ and ${\rm Im}(\rho_1) = {\rm pr}(\Gamma_1)$. So we want to prove
that $\overline{{\rm pr}(\Gamma_1)} = H_1$.

Let $h\in H_1$ be given. Since $H_1 \subset\overline{{\rm Im}(\rho)} =
\overline{{\rm pr}(\Gamma)}$, there exist a sequence
$(\gamma_n)_n\subset\Gamma$ such that $({\rm pr}(\gamma_n))_n$ converges
to $h$. But $H_1$ is open in $\overline{{\rm Im}(\rho)}$, because it is
its identity component. Hence, after dropping finitely many terms
from the sequence, we can assume that $({\rm pr}(\gamma_n))_n \subset
H_1$. This implies that $(\gamma_n)_n \subset G\times H_1$ and so
$(\gamma_n)_n \subset (G\times H_1)\cap\Gamma = \Gamma_1$. Hence,
we obtain $\overline{{\rm Im}(\rho_1)} = H_1$, thus completing the proof
of Theorem A.

\section{Proof of Theorem B}
\label{section-thmB}

To prove Theorem B, we assume that $G$ and $M$ satisfy its
hypotheses. Hence, we can choose $\widehat{M}$, $L=G\times H$,
$\Gamma$ and $K$ as in the conclusions of Theorem~A. First observe
that we can assume that $H \neq e$, since otherwise the
conclusions of Theorem B are trivially satisfied. In particular,
we will assume that the $G$-action on $M$ is not transitive.

We state the following result. We refer to
\cite{Zimmer-IHES-foliations} for the notion of compatible
transverse measure.

\begin{lemma}
\label{prop-ergodic} Every nontransitive locally free ergodic
action that preserves a compatible finite smooth measure
transverse to the orbits is properly ergodic.
\end{lemma}
\Proof The orbits define a foliation each of whose leaves
intersect a foliated chart in a countable number of plaques.
Hence, on a foliated chart every leaf intersects a transversal on
a countable and thus null set, since such transversals have
positive dimension. In particular, every orbit is null and so the
action is properly ergodic. \Endproof

By Lemma \ref{prop-ergodic} the $G$-action on $M$ is properly
ergodic. Then, the hypotheses of the main result in
\cite{Stuck-Zimmer} are satisfied and such result implies that the
$G$-action on $M$ is essentially free. Now we conclude from this
the following.

\begin{lemma}
\label{free-on-L-Gamma} The $G$-action on $\widehat{M}$ is
essentially free and the $G$-action on $L/\Gamma$ is free.
\end{lemma}
\Proof Since the $G$-action on $\widehat{M}$ is obtained as the
lift of the $G$-action on $M$ with respect to the covering map
$\widehat{M} \rightarrow M$, it follows that every $G$-orbit in
$\widehat{M}$ is a covering of a $G$-orbit in $M$. On the other
hand, since the $G$-action on $\widehat{M}$ is locally free, every
$G$-orbit on $\widehat{M}$ is a quotient of $G$ by a discrete
subgroup. If we consider such setup for a $G$-orbit $\OO$ in $M$
with trivial stabilizers and denote with $\widehat{\OO}$ a
$G$-orbit in $\widehat{M}$ that lifts $\OO$ we obtain a
commutative diagram of covering maps as follows:
$$
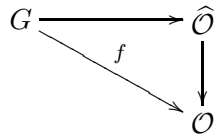
\xymatrix{ G \ar[rr] \ar[drr]^f  &    &  \widehat{\OO} \ar[d]   \\
&   &   \OO
 }
$$
where $f$ is a diffeomorphism.

Such diagram induces a corresponding commutative diagram of
fundamental groups given by:
$$
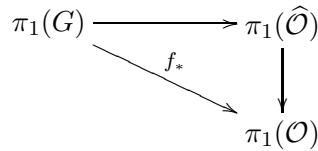
\xymatrix{ \pi_1(G) \ar[rr] \ar[drr]^{f_*} &
            &  \pi_1(\widehat{\OO}) \ar[d]   \\
&   &   \pi_1(\OO)
 }
$$
which can hold only if $\pi_1(\widehat{\OO}) = \pi_1(\OO)$, since
$f_*$ is an isomorphism and the vertical arrow is an inclusion.
Then the covering $\widehat{\OO} \rightarrow \OO$ is trivial,
which implies that the arrow $G \rightarrow \widehat{\OO}$ above
is a diffeomorphism and so the $G$-orbit $\widehat{\OO}$ has
trivial stabilizers.

Hence we have proved that every $G$-orbit in $M$ with trivial
stabilizers lifts to $\widehat{M}$ to $G$-orbits with trivial
stabilizers, and so the $G$-action on $\widehat{M}$ is essentially
free.

A similar argument proves that the $G$-action on $L/\Gamma$ is
essentially free. For this we use the fact that, as found in the
proof of Theorem A, $L/\Gamma$ is a\break $G$-invariant reduction of the
orthonormal frame bundle transverse to the foliation by
$G$-orbits, so that the $G$-orbits in $L/\Gamma$ are coverings of
their projections\break\vskip-12pt\noindent onto $\widehat{M}$, on which we just proved
essential freeness. Then we observe that for any given $l_1,l_2
\in L$, the stabilizers in $G$ of the $\Gamma$-classes of such
points satisfy $G_{l_2\Gamma} = l_2 l_1^{-1}G_{l_1\Gamma}l_1
l_2^{-1}$. Hence essential freeness of the $G$-action on
$L/\Gamma$ implies freeness for such action. \Endproof

Let $K_0$ be a maximal compact subgroup of $G$. By Lemma
\ref{free-on-L-Gamma} the\break $K_0$-action on $L/\Gamma$ is free and
so $K_0\backslash L/\Gamma$ is a compact connected manifold.
Moreover, the $G$-orbits in $L/\Gamma$ induce a foliation in
$K_0\backslash L/\Gamma$ that carries a leafwise Riemannian metric
so that each leaf is isometric to the simply connected Riemannian
symmetric space $K_0\backslash G$. The simply connectedness of
each leaf in $K_0\backslash L/\Gamma$ follows from Lemma
\ref{free-on-L-Gamma}. Also, by Theorem A, the $G$-action on
$L/\Gamma$ is topologically transitive. From this, it is easy to
check that the foliation in $K_0\backslash L/\Gamma$, by
Riemannian symmetric spaces, has a dense leaf. Finally, since $L =
G\times H$, the foliation on $K_0\backslash L/\Gamma$ carries a
transverse Lie structure modelled on $H$.

The previous discussion shows that the hypotheses in Theorem A in
\cite{Zimmer-Lie} are satisfied and such result implies that $H$
is semisimple, and so $L = G\times H$ is semisimple as well.

In what follows, for every connected semisimple Lie group $F$ we
denote with $F^{is}$ the minimal connected closed normal subgroup
of $F$ such that $F/F^{is}$ is compact. Then we clearly have
$G\subset L^{is}$.

On the other hand, by the proof of Theorem A, $\Gamma$ projects
densely into $H$ by the natural projection $L = G\times H
\rightarrow H$, which implies that $G\Gamma$ is dense in~$L$.
Hence, $L^{is}\Gamma$ is dense in $L$ as well. By Corollary~5.17
in \cite{Raghunathan} it follows that $\Gamma Z(L)$ is discrete
(see also Lemma 6.1 in page 329 in \cite{Margulis}). Hence,
$\Gamma Z(H)$ is discrete as well and, since it contains $\Gamma$,
it is also a lattice. This clearly implies that $\Gamma$ has
finite index in $\Gamma Z(H)$. But it is easy to prove that there
is a bijection:
$$
            (\Gamma Z(H))/\Gamma \cong Z(H)/(\Gamma\cap Z(H))
$$
and so $\Gamma\cap Z(H)$ has finite index in $Z(H)$.

Now denote $Z = \Gamma\cap Z(H)$, and observe that there is an
equivariant diffeomorphism:
$$
  \frac{L/Z}{\Gamma/Z} \cong L/\Gamma
$$
as well as isomorphisms:
$$
  KZ/Z \cong K/(K\cap Z) \cong K
$$
where the latter follows from $K\cap\Gamma = e$, which is a
consequence of the freeness of the $K$-action on $L/\Gamma$.

Hence, if we replace $L = G\times H$, $H$, $\Gamma$ and $K$ with
$L/Z = G\times H/Z$, $H/Z$, $\Gamma/Z$ and $KZ/Z$, then the
corresponding $\widehat{M}$  for the new choices is given by the
double coset $(KZ/Z)\backslash (L/Z)/(\Gamma/Z)$, which is easily
seen to be $G$-equivariantly diffeomorphic to our original choice
$K\backslash L/\Gamma$. Moreover, it is also clear that all the
properties proved so far still hold. Also, since we modded out by
$Z$, which has finite index in $Z(H)$, it follows that for our new
choices the group $L$ has finite center.

To complete to proof of Theorem B we will show that, for the above
choices, $\Gamma$ is irreducible and $L$ is isotypic. To achieve
this, we state the following general result.

\begin{lemma}
\label{lattice-irreducible} Let $F$ be a connected semisimple Lie
group with finite center and $\Lambda$ a lattice in $F$ such that
$F^{is}\Lambda$ is dense in $F$. Then the following conditions are
equivalent.
\begin{itemize}
    \item[{\rm (a)}] For every closed connected noncompact normal
    subgroup $N$ of $F${\rm ,} the group $N\Lambda$ is dense in $F$.
    \item[{\rm (b)}] For every closed connected noncompact normal proper
    subgroups $F'${\rm ,} $F''$ of $F$ that satisfy $F = F' F''$ with $F'\cap
    F''$ finite{\rm ,} the group $(\Lambda\cap F')(\Lambda\cap F'')$ has
    infinite index in $\Lambda$.
\end{itemize}
\end{lemma}
\Proof This result is essentially contained in Chapter V of
\cite{Raghunathan} although it is not explicitly stated. The proof
is obtained from the results in \cite{Raghunathan} as follows.

When $F = F^{is}$, the equivalence is part of Corollary 5.21 in
\cite{Raghunathan} and the definition of irreducible lattice in
this reference. The proof of such corollary ultimately depends on
Theorem 5.5 in \cite{Raghunathan}.

We then observe that our condition on $\Lambda$ implies that it
satisfies property (SS) (see \cite{Raghunathan} for the
definition) and so the hypotheses of Theorem 5.26 in
\cite{Raghunathan} are satisfied. In our case, it is easy to see
that, if we apply the arguments that prove Corollary 5.21 in
\cite{Raghunathan} using Theorem 5.26 in \cite{Raghunathan}
instead of Theorem 5.5 in \cite{Raghunathan}, we conclude the
equivalence between (a) and (b). \Endproof

The previous lemma shows that the definition of irreducible
lattice in a connected semisimple Lie group with finite center
(that may admit compact factors) as found in \cite{Morris}, which
is given by condition (a), is equivalent to condition (b). We now
use this to prove that our lattice $\Gamma$ is irreducible in $L$,
with irreducibility as defined in \cite{Morris} to be able to
apply results therein.

First recall that both $G\Gamma$ and $L^{is}\Gamma$ are dense in
$L$. By Lemma \ref{lattice-irreducible}, if $\Gamma$ is not
irreducible, then there exist closed connected noncompact normal
proper subgroups $L'$, $L''$ of $L$ that satisfy $L = L' L''$ with
$L'\cap L''$ finite and such that the group $(\Gamma\cap
L')(\Gamma\cap L'')$ has finite index in $\Gamma$. It follows easily
that $G(\Gamma\cap L')(\Gamma\cap L'')$ is dense in $L$ as well.
Since $G$ is a simple factor of $L$, it is either contained in $L'$
or in $L''$. If $G\subset L'$, then we have $G(\Gamma\cap
L')(\Gamma\cap L'') \subset L'(\Gamma\cap L'')$, where the latter is
a closed proper subgroup of $L$, which yields a contradiction. An
analogous contradiction is obtained if $G\subset L''$ and so
$\Gamma$ is irreducible.

Finally, since $L$ admits an irreducible lattice and is not
isogenous to either ${\rm SO}(1,n)\times C$ or ${\rm SU}(1,n)\times C$, for
any nontrivial connected compact group $C$, it follows from
Problem 6.8 in \cite{Morris} that $L$ is isotypic. This concludes
the proof of Theorem B.

\section{Proof of Theorem C}
\label{section-thmC}

The implication (2)$\Rightarrow$(1) is trivial by considering a
pseudoRiemannian metric on $K\backslash L/\Gamma$ as in
Proposition \ref{lattice-isotypic}.

Conversely, if we assume (1), then we can apply Theorem B to
conclude that $\widehat{X}$ has a finite covering of the form
$K\backslash L/\Gamma$ to which the $G$-action lifts. But then
$K\backslash L/\Gamma$ is a finite covering of $X$ and (2) holds.

\section{Further results and consequences}
\label{section-further}

A careful examination of the proofs in the previous sections show
that our assumption of compactness for $M$ in Theorem A is used
only to ensure completeness of the Riemannian or antiRiemannian
structure transverse to the $G$-orbits. Hence, the same proof
allows to obtain the following result.

\demo{\scshape Theorem A$'$}
{\it Let $G$ be a connected noncompact simple Lie group. If $G$ acts
faithfully and topologically transitively on a \/{\rm (}\/not necessarily
compact\/{\rm )}\/ manifold $M$ preserving a finite volume complete
pseudoRiemannian metric such that $n_0 = m_0${\rm ,} then the $G$-action
on $M$ is ergodic and engaging{\rm ,} and there exist\/{\rm :}\/
\begin{description}
    \item[{\rm (1)}] a finite covering $\widehat{M}\rightarrow M${\rm ,}
    \item[{\rm (2)}] a connected Lie group $L$ that contains $G$ as a factor{\rm ,}
    \item[{\rm (3)}] a discrete subgroup $\Gamma$ of $L$ and
    a compact subgroup $K$ of $C_L(G)${\rm ,}
\end{description}
for which the $G$-action on $M$ lifts to $\widehat{M}$ so that
$\widehat{M}$ is $G$-equivariantly diffeomorphic to $K\backslash
L/\Gamma$. Furthermore{\rm ,} there is an ergodic and engaging
$G$-invariant finite smooth measure on $L/\Gamma$.}
\Enddemo

On the other hand, the proof of our Theorem B makes use of
Theorem~A in \cite{Zimmer-Lie}. The latter requires compactness
only to ensure the existence of a\break $G$-invariant finite smooth
measure on a suitable manifold (see Proposition 2.5 in
\cite{Zimmer-Lie}). Besides that, no further use of the
compactness condition is made in \cite{Zimmer-Lie}. For our setup,
to have a corresponding version of Theorem~A in \cite{Zimmer-Lie}
for finite volume manifolds, we need to ensure that the $G$-action
on $L/\Gamma$ has an invariant finite smooth measure, which is a
conclusion in our Theorem A$'$. The proof of Theorem B also applies
Molino's work on Riemannian and Lie foliations, which are usually
stated for compact manifolds. However, if we assume completeness
of the manifold $M$, then the orthonormal transverse frame bundle
that appears in our arguments is transversely complete and
Molino's work still applies. Hence, the same arguments used to
prove Theorem B allows to obtain the following result.

\demo{\scshape Theorem B$'$}
{\it Let $G$ be a connected noncompact simple Lie group with\break finite
center and $\rank_\R(G) \geq 2$. If $G$ acts faithfully and
topologically transitively on a noncompact manifold $M$ preserving
a finite volume complete pseudo\-Riemannian metric such that $n_0 =
m_0${\rm ,} then there exist\/{\rm :}\/
\begin{itemize}
    \item[{\rm (1)}] a finite covering $\widehat{M}\rightarrow M${\rm ,}
    \item[{\rm (2)}] a connected isotypic semisimple Lie group $L$ without
    compact factors and with finite center that contains
    $G$ as a factor{\rm ,}
    \item[{\rm (3)}] an irreducible lattice $\Gamma$ of $L$ and
    a compact subgroup $K$ of $C_L(G)${\rm ,}
\end{itemize}
for which the $G$-action on $M$ lifts to $\widehat{M}$ so that
$\widehat{M}$ is $G$-equivariantly diffeo\-morphic to $K\backslash
L/\Gamma$. Hence{\rm ,} up to fibrations with compact fibers{\rm ,} $M$ is\break
$G$-equivariantly diffeomorphic to $K\backslash L/\Gamma$ and
$L/\Gamma$.}
\Enddemo

The additional property for $L$ in (2), its lack of compact factors,
follows from the following facts. Being irreducible, the lattice
$\Gamma$ is arithmetic in $L$ since the latter has finite center and
is not isogenous to either ${\rm SO}(1,n)\times C$ or ${\rm SU}(1,n)\times C$
for a compact group $C$. On the other hand, arithmetic lattices
which are not cocompact can only occur in semisimple Lie groups
without compact factors. Such claims follow from Theorem 6.21 and
Corollary 6.55 in \cite{Morris}.

From the above we obtain the following result.

\demo{\scshape Theorem C$'$}
{\it Let $G$ be a connected noncompact simple Lie group with finite
center and $\rank_\R(G) \geq 2$. Assume that $G$ acts faithfully
on a noncompact manifold $X$. Then the following conditions are
equivalent.
\begin{itemize}
    \item[{\rm (1)}] There is a finite covering $\widehat{X} \rightarrow
    X$ for which the $G$-action on $X$ lifts to a topologically
    transitive $G$-action on $\widehat{X}$ that preserves a finite
    volume complete pseudoRiemannian metric such that $n_0 = m_0$.
    \item[{\rm (2)}] There is a connected isotypic semisimple Lie group
    $L$ without compact factors and with finite center that contains
    $G$ as a factor{\rm ,} an irreducible lattice
    $\Gamma$ of $L$ and a compact subgroup $K$ of $C_L(G)$ such that
    $K\backslash L/\Gamma$ is a finite covering of $X$
    with $G$-equivariant covering map.
\end{itemize}}

We end this section with detailed proofs of Theorems
\ref{Riemannian-split}, \ref{gromov-rep2} and
\ref{virtual-arithmetic-quotient}.

\demo{Proof of Theorem {\rm \ref{Riemannian-split}}} We present the arguments
for compact $M$, but the proof is similar for $M$ complete with
finite volume. By Corollary \ref{lattice1},  $M$ has a finite covering of the form $(G\times
\widetilde{N})/\Gamma_0$.  Furthermore, from the proof of Theorem A it follows that the group
$\Gamma_0$ has a finite index subgroup $\Gamma_1$ contained in $R(G)\times \Iso(\widetilde{N})$.
In particular, $M_1 = (G\times \widetilde{N})/\Gamma_1$ is a finite
covering space of $M$. Hence it is enough to define a metric on
$M_1$ from one on $G\times \widetilde{N}$ which is the direct
product of a right invariant metric on $G$ and the metric on
$\widetilde{N}$. \Endproof

{\it Proof of Theorem} \ref{gromov-rep2}. Let $\widehat{M}$ and $\Gamma$
be as in the conclusions of Theorems B or B$'$. By (2) and (3) in
Theorems B and B$'$, $\Gamma$ admits a linear representation with the
required properties. Then the result is a consequence of the
following facts: 1) $\pi_1(\widehat{M})$ has finite index in
$\pi_1(M)$, and 2) there is a surjective homomorphism
$\pi_1(\widehat{M}) \rightarrow \Gamma$. \Endproof

{\it Proof of Theorem} \ref{virtual-arithmetic-quotient}. Choose
$\widehat{M} = K\backslash L/\Gamma$ a finite covering of $M$ as
provided by Theorems B or B$'$. By general results, it follows that
the $G$-action on $M$ has finite entropy even if $M$ is not compact
(see \cite{Lubotzky-Zimmer1} and \cite{Zimmer-book}).

On the other hand, if we project $L$, $K$ and $\Gamma$ using the
adjoint representation, then it is easy to obtain an arithmetic
$G$-space (as defined in \cite{Lubotzky-Zimmer1}), say
$K'\backslash L'/\Gamma'$, which admits a $G$-equivariant covering
of $\widehat{M}$. Then the result is a consequence of the
following two facts: 1) for our setup, $A(K'\backslash L'/\Gamma')
= K'\backslash L'/\Gamma'$, 2) if $Y_1\rightarrow Y_2$ is a
$G$-equivariant covering map with both $G$-actions ergodic and
with finite entropy, then $A(Y_1) = A(Y_2)$. Both facts follow
easily from the results and definitions in
\cite{Lubotzky-Zimmer1}. \hfill\qed

\end{document}